\documentclass[12pt]{article}
\pdfoutput=1
\usepackage[ascii]{inputenc}
\usepackage[T1]{fontenc}
\usepackage[english]{babel}
\usepackage{amsmath,amssymb,amsfonts,textcomp,epsfig}
\usepackage{hyperref}
\usepackage{color}
\usepackage{calc}
 \usepackage[all]{xy}

\def\R{\mathbb{R}}
\def\C{\mathbb C}
\def\Z{\mathbb Z}

\def\P{\mathbb P}
\title{The holonomy group at infinity of the Painlev\'e VI Equation}
\author{
Bassem Ben Hamed \\
 Institut Sup\'erieur des Sciences Appliqu\'ees et de Technologie de
  Gab\`es\\ D\'epartement de Math\'ematiques\\
Rue Amor Ben El Khatab, 6029 Gab\`es, Tunisie\and
Lubomir Gavrilov \\
Institut de Math\'{e}matiques de Toulouse, UMR 5219\\
  Universit\'{e}  de Toulouse,  31062 Toulouse,  France\and
Martine Klughertz \\
Institut de Math\'{e}matiques de Toulouse, UMR 5219\\
  Universit\'{e}  de Toulouse,  31062 Toulouse,  France}

\begin{document}
\maketitle
\newtheorem{definition}{Definition}
\newtheorem{remark}{Remark}
\newtheorem{theorem}{Theorem} 
\newtheorem{lemma}{Lemma}
\newtheorem{proposition}{Proposition}
\newtheorem{corollary}{Corollary}
\vspace{5mm} \noindent 2000 MSC scheme numbers: 70H07, 34M55, 37J30

\begin{abstract}
We prove that the holonomy group at infinity of the Painlev\'{e} VI equation is
virtually commutative.
\end{abstract}

\section{Introduction}
\label{section1} \noindent The sixth Painlev\'e equation
$\left(\mathbf{PVI}\right)$
\begin{eqnarray}\label{p6}
\nonumber \dfrac{d^2 \lambda}{dt^2} &= &\dfrac{1}{2}\left(\frac{1}{\lambda}+ \dfrac{1}
{\lambda - 1} + \dfrac{1}{\lambda - t}\right)\left(\dfrac{d\lambda}{dt}\right)^2 -
\left(\dfrac{1}{t} + \dfrac{1}{t - 1}
 + \dfrac{1}{\lambda - t}\right)\dfrac{d\lambda}{dt}\\
& & + \dfrac{\lambda (\lambda - 1)(\lambda - t)}
 {t^2 (t - 1)^2} \left[ \alpha - \beta \dfrac{t}{\lambda ^2} + \gamma \dfrac{t - 1}{(\lambda - 1)^2}
 +  (\dfrac{1}{2} - \delta) \dfrac{t(t - 1)}{(\lambda - t)^2}\right ]
\end{eqnarray}
is a family of differential equations parameterized by
$(\alpha,\beta,\gamma,\delta)\in \mathbb{C}^4$. \emph{The purpose of the
present paper is to show that the holonomy group of $\left(\mathbf{PVI}\right)$
at infinity is virtually commutative.} The precise meaning is as follows. It is
straightforward to check that (\ref{p6}) is equivalent to a non-autonomous
Hamiltonian system (the so called sixth Painlev\'e system)
\begin{eqnarray}\label{SH}\left\{\begin{array}{ccl}
\dfrac{d\lambda}{dt}&=&\dfrac{\partial H}{\partial \mu},\\
\dfrac{d\mu}{dt}&=&- \dfrac{\partial H}{\partial \lambda}
\end{array}\right.
\end{eqnarray}
where
\begin{eqnarray}\label{H}
\nonumber H&:=&\dfrac{1}{t(t-1)}\left[\lambda(\lambda-1)(\lambda-t)\mu^2+\left\{\kappa_0(\lambda-1)(\lambda-t)\right.\right.\\
&&\left.\left.+\kappa_1\lambda(\lambda-t)+(\kappa_t+1)\lambda(\lambda-1)\right\}\mu+\kappa(\lambda-t)\right]
\end{eqnarray}
and
\begin{equation}\label{relations}
\alpha=\dfrac{1}{2}\kappa_{\infty}^2,\;\beta=\dfrac{1}{2}\kappa_{0}^2,\;\gamma=\dfrac{1}{2}\kappa_{1}^2,
\;\delta=\dfrac{1}{2}\kappa_{t}^2 ,
\kappa=\dfrac{1}{2}\left[\left(\kappa_0+\kappa_1+\kappa_t+1\right)^2-\kappa_{\infty}^2\right].
\end{equation}
The phase space of the above system is
$$\{(\lambda,\mu,t)\in \mathbb{C}^3 :
t\neq 0,1\}$$
 which we
partially compactify  to $M=\mathbb{P}^1\times \mathbb{P}^1 \times
\{\mathbb{C}\setminus \{0,1\}\}$. It is immediately seen that the projective
lines
$$
\Gamma_c=\{ \mu=\infty, t= c\} \subset M, c\neq 0,1
$$
are leaves of the one-dimensional foliation  
induced by (\ref{SH}) on $M$. On each leaf $\Gamma_c$ the foliation has four
singular points defined by $\lambda =0,1,c,\infty $. Let $P\neq 0,1,c,\infty$
be a point on $\Gamma_c$ and consider a germ of a cross-section
$(\mathbb{C}^2,0)$ to $\Gamma_c$ at $P$. The holonomy group $G$ at infinity is
then the image of the holonomy representation
\begin{equation}\label{gc}
\pi_1(\Gamma_c\setminus\{0,1,t,\infty\},P) \rightarrow
Di\!f\!f(\mathbb{C}^2,0).
\end{equation}
It is defined up to a conjugation by a diffeomorphism, depending on the
germ of cross-section and the initial point $P$. Our main result is
\begin{theorem}
\label{main} The holonomy group at infinity of the sixth Painlev\'{e} equation is
virtually commutative.
\end{theorem}
 Recall that a group $G$ is said to be virtually commutative,
provided that there is a normal commutative subgroup $G^0\subset G$, such that
$G/G^0$ is finite. The isomorphism class of the holonomy group $G$ along the
leaf $\Gamma_c$ has in fact a canonical meaning. As we shall see in section
\ref{compactification}, the leaf $\Gamma_c$ coincides with the divisor $D_0(c)$
in the Okamoto compactification\cite{okamoto79}  of the phase space of $\mathbf{PVI}$, see fig.\ref{fig2}. In particular, the
holonomy group along $\Gamma_c$ is isomorphic to the holonomy group along the
Okamoto divisor $D_0(c)$. The remaining divisors shown on fig.\ref{fig2} are topological cylinders, the associated holonomy has therefore one generator and is commutative.

The proof of Theorem \ref{main} is based on Lemma \ref{mainlemma2} which claims
that the local holonomies near the singular points of the leaf $\Gamma_c$ are involutions, as well on the algebraic Lemma
\ref{mainlemma}. Lemma \ref{mainlemma2} and Lemma \ref{mainlemma} suggest that
Theorem \ref{main} is related to the fact that the vertical divisor shown on
fig.\ref{fig2} belongs to the  Kodaira list of degenerate elliptic curves.

Let $E_k$ be the $k$-th order variational equation along $\Gamma_c$ and $G_k$ the associated differential Galois
group. $E_k$ defines a connection on the Riemann sphere  $\Gamma_c= \P\setminus \{0,1,c,\infty\}$ with four regular singular points at the punctures $\{0,1,c,\infty\}$. 
The monodromy group of $E_k$ represents the $k$-th order jet of the holonomy group along $D_0(c)$. We describe these monodromy groups in the simplest cases $k=1,2$ in 
 section \ref{higher}. It follows, for instance, that the monodromy group of $E_1$ is isomorphic to a semi direct product $\Z^2\rtimes \Z_2$, while   $G_1= \C^2\rtimes \Z_2$. In particular $G_1$, as well $G_2$, are virtually commutative.
This is a particular case of a general fact. According to Theorem \ref{main} the monodromy group of $E_k$ is virtually commutative for all $k$.
As its Zarisky closure  is $G_k$, then we also have

\begin{theorem}
\label{main2}
For every $k$ the differential Galois group $G_k$ is virtually commutative.
\end{theorem}

The present paper was motivated by the study of the Liouville non-integrability of the 
$\mathbf{PVI}$ system through the Ziglin-Morales-Ramis-Simo theory of non-integrability \cite{mora10,ramis07,mora01,morales99}.
This theory asserts that
integrability in a Liouville sense along a particular solution $\Gamma_c$ implies that the variational equation $E_1$, as
well all higher order variational equations $E_k$ along this solution, have
virtually commutative differential Galois groups. Indeed, in such a way the "semi-local" non-integrability in a neighborhood of some particular solutions and parameter values of the $\mathbf{PVI}$ system has been  recently proved by Horozov and
Stoyanova \cite{horozov07,stoy09}, see also Morales-Ruiz \cite{morales06}. To prove the non-integrability for all parameters we  need, however, an explicitly known particular solution which exists for all parameter values. The only such appropriate solution is the vertical divisor $\Gamma_c=D_0(c)$, defined in Theorem
\ref{main}. The result of Theorem \ref{main2} shows that, contrary to what we expected,  one can not prove the absence of a first integral of the $\mathbf{PVI}$ equation, by making use of the Ziglin-Morales-Ramis-Simo theory. It is an open question, whether the $\mathbf{PVI}$ equation has a first integral, meromorphic along the divisor $\Gamma_c$ "at infinity". This question, but in a more general setting, has been raised in \cite[section 7]{ramis07}.

Non-integrability or transcendency of solutions  is one of the central subjects in the study of the $\mathbf{PVI}$ equation. 
The fact that its general solution can not be reduced to a solution of a first order differential equation has been claimed already by Painlev\'e, and proved more recently by Watanabe \cite{wata98} and others. A different approach to the transcendency, going back to Drach and Vessiot, is to interpret it as an irreducibility of  the Galois groupoid defined by Malgrange, see \cite{malgrange, casa07,casa09,casa12a}. The irreducibility of the $\mathbf{PVI}$ equation in the sense of Drach-Vessiot-Malgrange has been shown by Cantat and Loray \cite[Theorem 7.1]{calo09}. It follows from these results that the $\mathbf{PVI}$ equation does not allow an additional rational first integral. The relation between the irreducibility of the Galois groupoid of a Hamiltonian system and the differential Galois group along a given algebraic solution is studied recently by Casale \cite{casa12b}. In this context, our Theorem \ref{main} comes at a first sight as a surprise.  The solution $\Gamma_c$ which we use is however rather special : it is an irreducible component of the anti-canonical divisor of the space of initial conditions, and hence it is invariant  under the action of the Galois groupoid. This leads to special properties of the Galois groupoid along $\Gamma_c$ too.

The paper is organized as follows. In section \ref{compactification}, we resume
briefly the Okamoto compactification of the phase space of $\mathbf{PVI}$ equation
\cite{okamoto79}. In section \ref{higher} we describe 
the monodromy group of the first and the second variational equation along $\Gamma_c$, in terms of complete elliptic integrals of first and second kind. These groups provide an approximation of the holonomy group along $\Gamma_c$.  
Our main
result, Theorem \ref{main}, is proved in section \ref{mainsection}.

\section{The Okamoto compactification}
\label{compactification}

Let $(E,\pi,B)$ be a complex-analytic fibration with base $B$, total space $E$
and projection $\pi : E\rightarrow B$. Consider a foliation $\mathcal{F}$ on
$E$ of dimension equal to the dimension of $B$. Following \cite{okamoto79} we
say that $\mathcal{F}$ is P-uniform, if for every leaf $\Gamma\subset E$ the
induced map
$$
\pi : \Gamma \rightarrow B
$$
is an analytic covering. Thus, for every initial point $e\in E$, and every
continuous path $\gamma\subset B$ starting at $b=\pi(e)$, there is a unique
continuous path $\tilde{\gamma}\subset E$ starting at $e$, which is a lift of
$\gamma$ with respect to $\pi$ (the "Painlev\'{e} property" of the foliation). The
analyticity of $\pi$ implies moreover that at each point $e\in E$ the leaf of
the foliation is transversal to the corresponding fiber of the fibration.

From now on we put
$$E= \{ (\lambda,\mu,t)\in \mathbb{C}^3 : t\neq 0,1\}, B =
\mathbb{C}\setminus \{0,1\}
$$
$$\pi : E \rightarrow B : (\lambda,\mu,t) \mapsto t
$$ being the natural projection. The system (\ref{SH}) defines a one-dimensional
foliation $\cal{F}$ on the total space $E$ which is not P-uniform, but can be
completed to a P-uniform foliation after an appropriate partial
compactification $\bar{E}$ of $E$.

\begin{figure}
\begin{center}
\includegraphics[width=12cm]{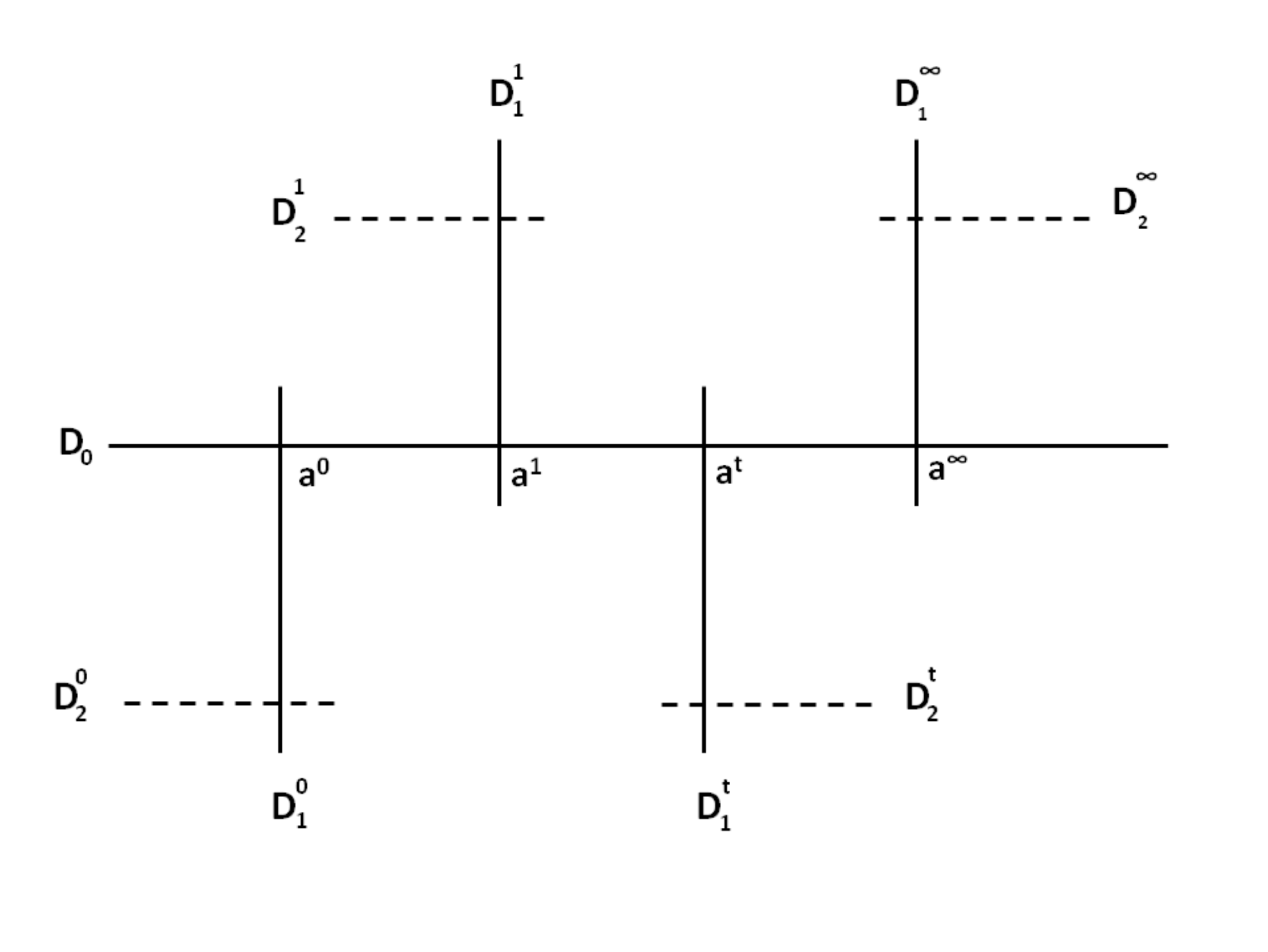}
\end{center}
\caption{The divisor $\bar{E}_t\setminus E_t$} \label{fig2}
\end{figure}

The main result of \cite{okamoto79} may be formulated as follows.
\begin{theorem}
\label{mainokamoto}
 There exists a canonical compact complex-analytic fibration
$(\bar{E},\bar{\pi},B)$, such that
\begin{itemize}
    \item $E\subset \bar{E}$, $\bar{\pi}|_E=\pi$
    \item Each fiber $\bar{E}_t= \bar{\pi}^{-1}(t)$ is compact
    \item $\bar{E}_t\setminus E_t$ is a union of nine transversal projective lines,
    as it is shown on fig.\ref{fig2}. The intersection points of the lines depend
    analytically on $t$.
    \item Let $D_t$ be the union of five solid lines shown on fig.\ref{fig2}.The
    foliation induced by (\ref{SH}) on $\tilde{E}=
    \bar{E}\setminus \cup_{t\in B} D_t$ is P-uniform with respect to the induced
    projection.
\end{itemize}
\end{theorem}
A similar result holds true for the remaining Painlev\'{e} equations
\cite{okamoto79}.\\
 \textbf{Remark.} $\bar{E}_t\setminus  D_t$ is the so called
"space of initial conditions" of the Painlev\'{e} VI equation which we describe
next.\\
\textbf{Sketch of the proof of Theorem \ref{mainokamoto}.}
Following \cite[Okamoto]{okamoto79}, define first the Hirzebruch surface
$\Sigma^{(2)}_{(\epsilon)}$, $\varepsilon \in \C$, using four charts $W_i= \C^2$,
 with local coordinates $(\lambda_i,\mu_i)$, $i=1,\ldots,4$, where
\begin{eqnarray}\left\{\begin{array}{ccc}\label{cartes}
\lambda_2=\lambda_1,\;\mu_2=\dfrac{1}{\mu_1} &\;\text{ in }\;&W_1\cap W_2,\\
\lambda_3=\dfrac{1}{\lambda_1},\;\mu_3=\epsilon\lambda_1-\lambda_1^2\mu_1 &\;\text{ in }\;&W_1\cap W_3,\\
\lambda_4=\lambda_3,\;\mu_4=\dfrac{1}{\mu_3} &\;\text{ in }\;&W_3\cap
W_4.\end{array}\right.
\end{eqnarray}
If $\epsilon\neq 0$, then the Hirzebruch surface $\Sigma^{(2)}_{(\epsilon)}$ is isomorphic to $\P^1\times\P^1$; otherwise it is isomorphic to
the tangent projective bundle of $\P^1$ with projection
\begin{eqnarray*}
\Sigma^{(2)}_{(\epsilon)}&\rightarrow &\P^1 \\
(\lambda_i,\mu_i)& \mapsto & \lambda_i  \; .
\end{eqnarray*}
The vector field (\ref{SH})  extends on the total space of the trivial bundle
\begin{equation}
\label{bundle} \Sigma^{(2)}_{(\epsilon)}\times B \stackrel{\pi}{\rightarrow} B,
B=\P^1\setminus\{0,1,\infty\}
\end{equation}
where $\epsilon=-\left(\kappa_0+\kappa_1+\kappa_t+\kappa_{\infty}+1\right)$. For instance, in  the chart $W_2$  it takes the form

\begin{eqnarray}\label{2}
 \left\{\begin{array}{ccl}
\mu_2\lambda_2^{\prime}&=&\dfrac{1}{t(t-1)}\left[2E(t,\lambda_2)+F(t,\lambda_2)\mu_2\right],\\
\mu_2^{\prime}&=&\dfrac{1}{t(t-1)}\left[E_{\lambda}(t,\lambda_2)+F_{\lambda}(t,\lambda_2)\mu_2+G\mu_2^2\right],
\end{array}\right.
\end{eqnarray}
where
\begin{eqnarray}\label{A}
 \left\{\begin{array}{l}
E(t,\lambda)=\lambda(\lambda-1)(\lambda-t),\\
F(t,\lambda)=\kappa_0(\lambda-1)(\lambda-t)+\kappa_1\lambda(\lambda-t)+(\kappa_t+1)\lambda(\lambda-1),\\
E_{\lambda}=\dfrac{\partial E}{\partial\lambda},\;F_{\lambda}=\dfrac{\partial F}{\partial\lambda},\\
G=-\dfrac{1}{2}\epsilon\left(\kappa_0+\kappa_1+\kappa_t-\kappa_{\infty}+1\right)=\kappa,
\end{array}\right.
\end{eqnarray}
The above meromorphic vector field 
induces a singular foliation on $\Sigma^{(2)}_{(\epsilon)}$ 
 having four one-parameter families of singular points
$S^{\theta}$, $\theta=0,1,\infty,t$  defined by
\begin{equation*}
S^{\theta}\cap \pi^{-1}(t)=a^{\theta}(t),
\end{equation*}
\begin{eqnarray*}
a^0(t)&=&\left\{(\lambda_2,\mu_2)=(0,0)\right\},\\
a^1(t)&=&\left\{(\lambda_2,\mu_2)=(1,0)\;\text{or}\;(\lambda_4,\mu_4)=(1,0)\right\},\\
a^t(t)&=&\left\{(\lambda_2,\mu_2)=(t,0)\;\text{or}\;(\lambda_4,\mu_4)=(\dfrac{1}{t},0)\right\},\\
a^{\infty}(t)&=&\left\{(\lambda_4,\mu_4)=(0,0)\right\}.
\end{eqnarray*}
Replace the Hirzebruch surface $\Sigma^{(2)}_{(\epsilon)}$ by
$\Sigma^{(2)}_{(\epsilon)}$ blown up at $a^{\theta}(t)$ for every $t\in B$.
This replaces each $a^{\theta}(t)$ by a projective line denoted
$D_1^\theta(t)$. The induced foliation has still four one-parameter families of
singular points which belong to $D_1^\theta(t)$. We blow up once again the
surfaces at these singular points to obtain the fibers $\bar{E}_t$ of the
fibration described in Theorem \ref{mainokamoto}, see fig.\ref{fig2}. The
remaining claims of the
Theorem follow by computation.$\Box$\\

\section{Higher order variational equations and their monodromy groups}
\label{higher}
 In this section, we consider the foliation $\mathcal{F}$ defined
by the vector field (\ref{SH}) on the total space of the fibration $(\bar{E},\bar{\pi},B)$, see Theorem \ref{mainokamoto}.
This foliation has in each fiber $\pi^{-1}(t)$ a vertical leaf $D_0(t)$, which in the   chart $W_2$ takes the form\begin{equation*}
D_0(t):\mu_2=0 .
\end{equation*}
According to (\ref{2}) the foliation $\mathcal{F}$ in the local chart $W_2$ is defined by
\begin{eqnarray}\label{2prime}
\left\{\begin{array}{ccl}
d\mu&=&\dfrac{\left[E_{\lambda}(t,\lambda)+F_{\lambda}(t,\lambda)\mu+G\mu^2\right]\mu}{2E(t,\lambda)+F(t,\lambda)\mu} d\lambda\\
dt&=&\dfrac{t(t-1)\mu}{2E(t,\lambda)+F(t,\lambda)\mu} d\lambda
\end{array}\right.
\end{eqnarray}
where $E$, $F$ and $G$ are given by (\ref{A}). Here, as well until the end of the paper, we replace for simplicity $\mu_2, \lambda_2$ by $\mu, \lambda$.\\

In this section we compute the first and the second variational equations of (\ref{2prime}) along $D_0(c)$ and study the corresponding monodromy groups. For this purpose we put, following \cite{horozov07,ramis07},
$$
t=c+ \varepsilon \eta_1 + \dfrac {\varepsilon^2} 2 \eta_2 + \dots, \quad \mu= \varepsilon\xi_1 + \dfrac{\varepsilon^2}2 \xi_2 + \dots, \quad \varepsilon \sim 0
$$
where $\eta_k = \eta_k(\lambda) ,\xi_k = \xi_k(\lambda)$ are unknown functions,
and substitute these expressions in (\ref{2prime}). Equating the coefficients of $\varepsilon^k$ we get a recursive system of  linear non-homogeneus equations on $(\eta_k,\xi_k)$ - the higher order variational equations. We note that these equations, except in the case $k=1$, are non-linear. In order to obtain a linear system we add suitable monomials in $\eta_i, \xi_j$, e.g. \cite{ramis07, mora10}. The fundamental matrices of solutions of these equations are then explicitly computed by the Picard method in terms of iterated integrals. This implies also a description of the corresponding monodromy matrices. In the next two subsections we carry out this procedure in the particular case of the first and the second variational equation.

\subsection{The first variational equation}
The first variational equation $E_1$ along $D_0$ is the linear
system
\begin{eqnarray}\label{ve1}
\left(\begin{array}{c}
\dot{\eta}_1\\
\dot{\xi}_1\end{array}\right)=\left(\begin{array}{cc}
0 & b(\lambda)\\
0 & a(\lambda)\end{array}\right)\left(\begin{array}{c}
\eta_1\\
\xi_1\end{array}\right),
\end{eqnarray}
where
\begin{eqnarray*}
a(\lambda)&=&\dfrac{E_{\lambda}(c,\lambda)}{2E(c,\lambda)},\\
b(\lambda)&=&\dfrac{c(c-1)}{2E(c,\lambda)}.
\end{eqnarray*}
The general solution of the system (\ref{ve1}) is  given by

\begin{eqnarray*}
\eta_1(\lambda)&=&c_1\int^\lambda_p \dfrac{c(c-1)d \lambda}{2\sqrt{\lambda(\lambda-1)(\lambda-c)}} +c_2,\\
\xi_1(\lambda)&=&c_1 \sqrt{\lambda(\lambda-1)(\lambda-c)}.
\end{eqnarray*}
where $(c_1,c_2)\in\C^2$ and $p\in \C$ is a fixed initial point.
The fundamental matrix of solutions
 \begin{eqnarray*} X(\lambda) =
\left(\begin{array}{cc}
\int^\lambda_p \dfrac{c(c-1)d \lambda}{2\sqrt{\lambda(\lambda-1)(\lambda-c)}} & 1\\
&  \\
\sqrt{\lambda(\lambda-1)(\lambda-c)} & 0\\
\end{array}\right)
\end{eqnarray*}
is multivalued, and the result of the analytic continuation of $X(.)$ along small loops making one turn around $\lambda= 0, 1, c$ respectively is
$$
X \rightarrow X T_0, X \rightarrow X T_1, X \rightarrow X T_c.
$$
The matrices $T_0, T_1, T_c$ generate the monodromy group of (\ref{ve1}) and can be computed as follows.
Let $S_c$ be the
 compact elliptic Riemann surface of the algebraic function $\sqrt{\lambda(\lambda-1)(\lambda-c)}$. It has 
 an  affine equation
\begin{equation}
\label{ec}
\{ (\lambda,y): y^2= \lambda(\lambda-1)(\lambda-c) \} . 
\end{equation}
The one-form $$\dfrac{d \lambda}{\sqrt{\lambda(\lambda-1)(\lambda-c)}}$$ is holomorphic on $S_c$
and hence
$X(.)$ can be seen as a globally multivalued, but locally meromorphic matrix function on $S_c$. This implies that
$$
T_0^2=T_1^2=T_c^2= \left(\begin{array}{cc}
1 & 0\\
0 & 1\\
\end{array}\right)
$$
and hence
\begin{eqnarray}
\label{ttt}
T_0= \left(\begin{array}{cc}
-1 & 0\\
\alpha_0 & 1\\
\end{array}\right), T_1= \left(\begin{array}{cc}
-1 & 0\\
\alpha_1 & 1\\
\end{array}\right), T_c= \left(\begin{array}{cc}
-1 & 0\\
\alpha_c & 1\\
\end{array}\right).
\end{eqnarray}
The constants $\alpha_0, \alpha_1, \alpha_c$  depend on  the initial point $p$ and can be determined as follows. The matrix 
$$
T_0 T_1 = \left(\begin{array}{cc}
1 & 0\\
\alpha_1-\alpha_0 & 1\\
\end{array}\right)
$$ represents the monodromy of $X(.)$ along a closed loop on the $\lambda$-plane, which lifts,  on on the Riemann surface of $\sqrt{\lambda(\lambda-1)(\lambda-c)}$ to a closed loop too, which we denote
$\gamma$.
The monodromy of the fundamental matrix $X$ along this loop is $T_0 T_1$ and we have 
\begin{eqnarray*}X \rightarrow X T_0 T_1= X+ 
\left( \begin{array}{cc}\Pi& 0  \\
0 & 0\end{array} \right)
\end{eqnarray*}
where
$$
\Pi= \int_\gamma  \dfrac{c(c-1)d \lambda}{2\sqrt{\lambda(\lambda-1)(\lambda-c)}}  =  \int_0^1\dfrac{c(c-1)d \lambda}{\sqrt{\lambda(\lambda-1)(\lambda-c)}} 
$$
is a period of the holomorphic one-form on $S_c$. Therefore 
$$
\alpha_1-\alpha_0=  \int_0^1\dfrac{c(c-1)d \lambda}{\sqrt{\lambda(\lambda-1)(\lambda-c)}} 
$$
and in a similar way
$$
\alpha_c -\alpha_0=  \int_0^c \dfrac{c(c-1)d \lambda}{\sqrt{\lambda(\lambda-1)(\lambda-c)}} .
$$
Finally, taking the limit $p \rightarrow 0$ we obtain $\alpha_0=0$. The monodromy group of (\ref{ve1}) is therefore 
\begin{equation}
\label{mon1}
<T_0,T_1,T_c> =  \{ \left( \begin{array}{cc}\pm 1& 0  \\
p\Pi_1+q\Pi_2 &1 \end{array} \right) : p,q\in \Z \}
\end{equation}
where
$$
\Pi_1 = \int_0^1\dfrac{c(c-1)d \lambda}{\sqrt{\lambda(\lambda-1)(\lambda-c)}} , \quad \Pi_2 = \int_0^c\dfrac{c(c-1)d \lambda}{\sqrt{\lambda(\lambda-1)(\lambda-c)}} 
$$
are the fundamental periods of the elliptic surface $S_c$. As $\Pi_1, \Pi_2$ are linearly independent over $\R$, then we obtain
\begin{proposition} For every $c\neq 0,1,t$,
the monodromy group (\ref{mon1}) of the first variational equation (\ref{ve1})   is  isomorphic to the semidirect product $\Z^2 \rtimes \Z_2$, where $\Z_2=\Z/2\Z$.
\end{proposition}
It is well known that for a Fuchs type equation
the Zariski closure of the monodromy group is the differential Galois group (e.g. \cite{tret79}). Therefore the Galois group 
$G_1$ of (\ref{ve1}) is 
$$
G_1= \C \rtimes \Z_2 = \{ \left( \begin{array}{cc}\pm 1& 0  \\
z&1 \end{array} \right) : z\in \C \}.
$$
In particular,  $G_1$ is virtually commutative. 

\subsection{The second variational equation}

The second variational equation $E_2$ along the divisor $D_0(c)$
reads
\begin{eqnarray}\label{ve2}
\nonumber \dot{\xi}_2&=&d(\lambda)\left(\xi_1\right)^2+e(\lambda)\xi_1\eta_1+a(\lambda)\xi_2,\\
\dot{\eta}_2&=&f(\lambda)\left(\xi_1\right)^2+g(\lambda)\xi_1\eta_1+b(\lambda)\xi_2,
\end{eqnarray}
where
\begin{eqnarray*}
a(\lambda)&=&\dfrac{E_{\lambda}(c,\lambda)}{2E(c,\lambda)},\\
b(\lambda)&=&\dfrac{c(c-1)}{2E(c,\lambda)},\\
d(\lambda)&=&\dfrac{2E(c,\lambda)F_{\lambda}(c,\lambda)-E_{\lambda}(c,\lambda)F(c,\lambda)}{4(E(c,\lambda))^2},\\
e(\lambda)&=&\dfrac{-(2\lambda-1)E(c,\lambda)+\lambda(\lambda-1)E_{\lambda}(c,\lambda)}{2(E(c,\lambda))^2},\\
f(\lambda)&=& -\dfrac{c(c-1)F(c,\lambda)}{4(E(c,\lambda))^2},\\
g(\lambda)&=&\dfrac{(2c-1)E(c,\lambda)+c(c-1)\lambda(\lambda-1)}{2(E(c,\lambda))^2}.
\end{eqnarray*}
Having computed $\mu_1, \xi_1$, this is  a linear non-homogeneous equation in $\mu_2, \xi_2$, but it is also 
equivalent to the linear  system:
\begin{eqnarray}\label{linear}
\left(\begin{array}{c}
\dot{\eta}_2\\
\dot{\xi}_2\\
\dot{u}_1\\
\dot{v}_1
\end{array}\right)=\left(\begin{array}{cccc}
0 & b(\lambda) & g(\lambda) & f(\lambda)\\
0 & a(\lambda) & e(\lambda) & d(\lambda)\\
0 & 0 & a(\lambda) & b(\lambda)\\
0 & 0 & 0 & 2a(\lambda)
\end{array}\right)\left(\begin{array}{c}
\eta_2\\
\xi_2\\
u_1\\
v_1
\end{array}\right),
\end{eqnarray}
where $u_1=\xi_1\eta_1$ et $v_1=\left(\xi_1\right)^2$.\\
The substitution
\begin{equation*}
\sigma_2=\dfrac{\xi_2}{\sqrt{\lambda(\lambda-1)(\lambda-c)}},\;u_2=\dfrac{u_1}{\sqrt{\lambda(\lambda-1)(\lambda-c)}},\;v_2=\dfrac{v_1}{\lambda(\lambda-1)(\lambda-c)},
\end{equation*}
transforms (\ref{linear}) to a strictly upper triangular form
\begin{eqnarray}\label{linear2}
\left(\begin{array}{c}
\dot{\eta}_2\\
\dot{\sigma}_2\\
\dot{u}_2\\
\dot{v}_2
\end{array}\right)=\left(\begin{array}{cccc}
0 & a_{12}(\lambda) & a_{13}(\lambda) & a_{14}(\lambda)\\
0 & 0 & a_{23}(\lambda) & a_{24}(\lambda)\\
0 & 0 & 0 & a_{34}(\lambda)\\
0 & 0 & 0 & 0
\end{array}\right)\left(\begin{array}{c}
\eta_2\\
\sigma_2\\
u_2\\
v_2
\end{array}\right),
\end{eqnarray}
where
\begin{eqnarray*}
a_{12}(\lambda)&=&\sqrt{\lambda(\lambda-1)(\lambda-c)}\;b(\lambda),\\
a_{13}(\lambda)&=&\sqrt{\lambda(\lambda-1)(\lambda-c)}\;g(\lambda),\\
a_{14}(\lambda)&=&\lambda(\lambda-1)(\lambda-c)\;f(\lambda),\\
a_{23}(\lambda)&=&e(\lambda),\\
a_{24}(\lambda)&=&\sqrt{\lambda(\lambda-1)(\lambda-c)}\;d(\lambda)\\
a_{34}(\lambda)&=& a_{12}(\lambda) .
\end{eqnarray*}
The linear system (\ref{linear2}) is solved recursively in terms of iterated integrals. Namely, for differential forms  $\omega_i(x)=f_i(x)dx$ on the interval $[0,1]$ define the linear iterated integrals
$$
\int_0^1 \omega_1 \omega_2 = \int_0^1 f_1(y) (\int_0^y f_2(x)dx) dy,
$$
$$
\int_0^1 \omega_1 \omega_2 \omega_3= \int_0^1 f_1(z)  [\int_0^z f_2(y) (\int_0^y f_3(x) dx)dy] dz .
$$
Integrals of higher order and along a path on a  Riemann surface are defined in a similar way, e.g. \cite{hain87}. The fundamental matrix of solutions of the linear system (\ref{linear2}) takes the form
\begin{eqnarray*}
Y(\lambda)= I+ \int_P^\lambda J + \int_P^\lambda J^2+ \int_P^\lambda J^3 =
\left(\begin{array}{cccc}
1 & Y_{12}(\lambda) & Y_{13}(\lambda) & Y_{14}(\lambda)\\
0 & 1 & Y_{23}(\lambda) & Y_{24}(\lambda)\\
0 &0& 1 & Y_{34}(\lambda)\\
0 & 0 & 0 & 1
\end{array}\right),
\end{eqnarray*}
where
\begin{eqnarray*}
J = 
\left(\begin{array}{cccc}
0 & \omega_{12}(\lambda) & \omega_{13}(\lambda) & \omega_{14}(\lambda)\\
0 & 0 & \omega_{23}(\lambda) & \omega_{24}(\lambda)\\
0 &0& 0 & \omega_{34}(\lambda)\\
0 & 0 & 0 & 0
\end{array}\right), \quad \omega_{ij}= a_{ij}(\lambda) d \lambda .
\end{eqnarray*}
and hence
\begin{eqnarray*}
Y_{12}(\lambda)&=&\int_P^\lambda \omega_{12},\\
Y_{13}(\lambda)&=&\int_P^\lambda\omega_{12}\omega_{23}+\int_P^\lambda\omega_{13},\\
Y_{14}(\lambda)&=&\int_P^\lambda\omega_{12}\omega_{23}\omega_{34}+\int_P^\lambda\omega_{12}\omega_{24}+\int_P^\lambda\omega_{13}\omega_{34}+\int_P^\lambda\omega_{14},\\
Y_{23}(\lambda)&=&\int_P^\lambda\omega_{23},\\
Y_{24}(\lambda)&=&\int_P^\lambda\omega_{23}\omega_{34}+\int_P^\lambda\omega_{24},\\
Y_{34}(\lambda)&=&\int_P^\lambda\omega_{34}.
\end{eqnarray*}

As for the first variational equation, the fundamental matrix $Y(\lambda)$ is a globally multivalued, but locally meromorphic matrix function on the elliptic curve $S_c$, (\ref{ec}) (after removing eventually the points $\lambda=0,1,c, \infty$). The monodromy matrix $T_\alpha$ of $Y$ along a  closed path $\alpha$  
 is given therefore by the same matrix $Y$, in which the integrals $\int_P^\lambda$ are replaced by $\int_{\alpha}$. A more careful analysis will show, however, that when $\alpha$ can be lifted to a closed loop on $S_c$, then the double and triple iterated integrals in $T_\alpha$ are reduced to usual complete elliptic  integrals of first and second kind. This would imply the involutivity of the monodromy operators $T_0, T_1, T_c$ (defined in the receding sub-section) as well the virtual commutativity of the monodromy group of $E_2$.

Let $\alpha$ be a closed path on $D_0(c)$ which lifts to a closed path on the elliptic curve $S_c$, where 
\begin{eqnarray*}
S_c&\rightarrow &D_0(c)\\
(\lambda,\mu)& \mapsto & \lambda
\end{eqnarray*}
 is a double ramified covering over $\lambda = 0,1,c, \infty$. Denote the monodromy matrix of the second variational equation $(\ref{linear})$ along $\alpha$ by $T_\alpha$. 
\begin{proposition} 
\label{pe2} 
The entries of the monodromy matrix $T_\alpha$ along a closed loop $\alpha$ on the elliptic surface $S_c$ are quadratic polynomials in the complete elliptic integrals of first and second kind along $\alpha$.
\end{proposition}
Until the end of this subsection we sketch the proof  the Proposition \ref{pe2}.
Note that $\omega_{23}=d\,h $ is an exact form,
where $$h(\lambda)=-\dfrac{\lambda(\lambda-1)}{2E(c,\lambda)}= -\dfrac 1{2(\lambda-c)} .
 $$
 This combined with the identity
$$
\int_P^\lambda \omega_1 \omega_2 + \int_P^\lambda \omega_2 \omega_1 = \int_P^\lambda \omega_1 \int_P^\lambda\omega_2
$$
allows to express the iterated integrals of length two and three via usual Riemann integrals of meromorphic one-forms. Indeed,
 for
every differential 1-form $\omega$, we have
\begin{equation}
\label{double}
\int_{P}^\lambda\omega_{23}\omega=h(\lambda)\int_{P}^\lambda\omega-\int_{P}^\lambda h(\lambda)\omega .
\end{equation}
It follows that $Y_{13}(\lambda)$ and $ Y_{24}(\lambda)$ are in fact Riemann integrals along meromorphic differential forms. The reader may check that these one-forms have no residues on $S_c$. It remains to analyze $ Y_{14}(\lambda)$ .
Using  (\ref{double}), we obtain
\begin{equation}
\label{triple}
\int_{P}^\lambda\omega_{12}\omega_{23}\omega_{12}=\int_{P}^\lambda\left(\omega_{12}h(\lambda)\right)
\omega_{12}-\int_{P}^\lambda\omega_{12}\left(h(\lambda)\omega_{12}\right).
\end{equation}
which implies
\begin{eqnarray*}
\int_{P}^\lambda \omega_{12}\omega_{23}\omega_{12}+\int_{P}^\lambda\omega_{12}\omega_{24}+\int_{P}^\lambda\omega_{13}\omega_{12}&=&\int_{P}^\lambda\omega_{12}\overline{\omega}\\
& &+( \int_{P}^\lambda\omega_{12} )(\int_{P}^\lambda h \; \omega_{12}+\int_{P}^\lambda\omega_{13}),
\end{eqnarray*}
where
$\overline{\omega}:=\left(-2h(\lambda)b(\lambda)+d(\lambda)-g(\lambda)\right)\sqrt{E(c,\lambda)}\;d\lambda$.
The two forms $\omega_{12}$ and $\overline{\omega}$ have dependent cohomology
classes in $H^{1}_{DR}\left(\Gamma\right)$, where $S_c$ is the elliptic
curve(\ref{ec}). In fact,
\begin{eqnarray*}
\omega_{12}=\dfrac{c(c-1) d\lambda}{2 \sqrt{E(c,\lambda)}}\;\text{ and
}\;\overline{\omega}=\left(\dfrac{1-2c}{2}\right)\dfrac{d\lambda}{\sqrt{E(c,\lambda)}}+
\dfrac{1}{2}d\left(\dfrac{F(c,\lambda)}{\sqrt{E(c,\lambda)}}\right) .
\end{eqnarray*}
Thus everything is reduced to quadratic expressions in suitable Riemann integrals along meromorphic differential form without residues on $S_c$.
(the latter claim is straightforward to check). From this the the Proposition follows. $\Box$
\begin{corollary}
Let $\alpha, \beta$ be closed loops on $S_c$ starting at the same point, so they can be composed. Then $T_\alpha T_\beta = T_\beta T_\alpha$ and the monodromy group of $E_2$ is virtually commutative.
\end{corollary}
Indeed, as the homology class of the loop $\alpha \beta \alpha^{-1} \beta^{-1} $ is zero, then 
$$T_\alpha T_\beta T_{\alpha}^{-1} T_{\beta}^{-1}= T_{\alpha \beta \alpha^{-1} \beta^{-1} }$$ is the unit matrix.
Note that the  monodromy operators along closed loops on $S_c$  generate a subgroup of the monodromy group of the second variational equation $(\ref{linear})$ of finite index. Therefore, as expected, the monodromy group is virtually commutative.

\section{The holonomy at infinity}
\label{mainsection}
 Let $\Gamma$ be a leaf of the Painlev\'{e} VI foliation and
$(\mathbb{C}^2,0)$ be a germ of a cross-section to $\Gamma$ at some regular
point $P$. Each homotopy class of closed loops $\gamma\subset\Gamma$, starting at
$P$, defines a germ of a diffeomorphism
$$
h_\gamma : \mathbb{C}^2,0 \rightarrow \mathbb{C}^2,0
$$
and a homomorphism (the holonomy representation of the fundamental group of $\Gamma$)
$$
\pi_1(\Gamma, P) \rightarrow Di\!f\!\!f(\mathbb{C}^2,0): \gamma \mapsto
h_\gamma .
$$
The holonomy group of the foliation along $\Gamma$ is the image of this map
(which will be confounded with the representation itself ). Different points in
the leaf and different cross-sections give rise to representations conjugated
by germs of holomorphic diffeomorphisms.
\begin{figure}
\begin{center}
\includegraphics[width=12cm]{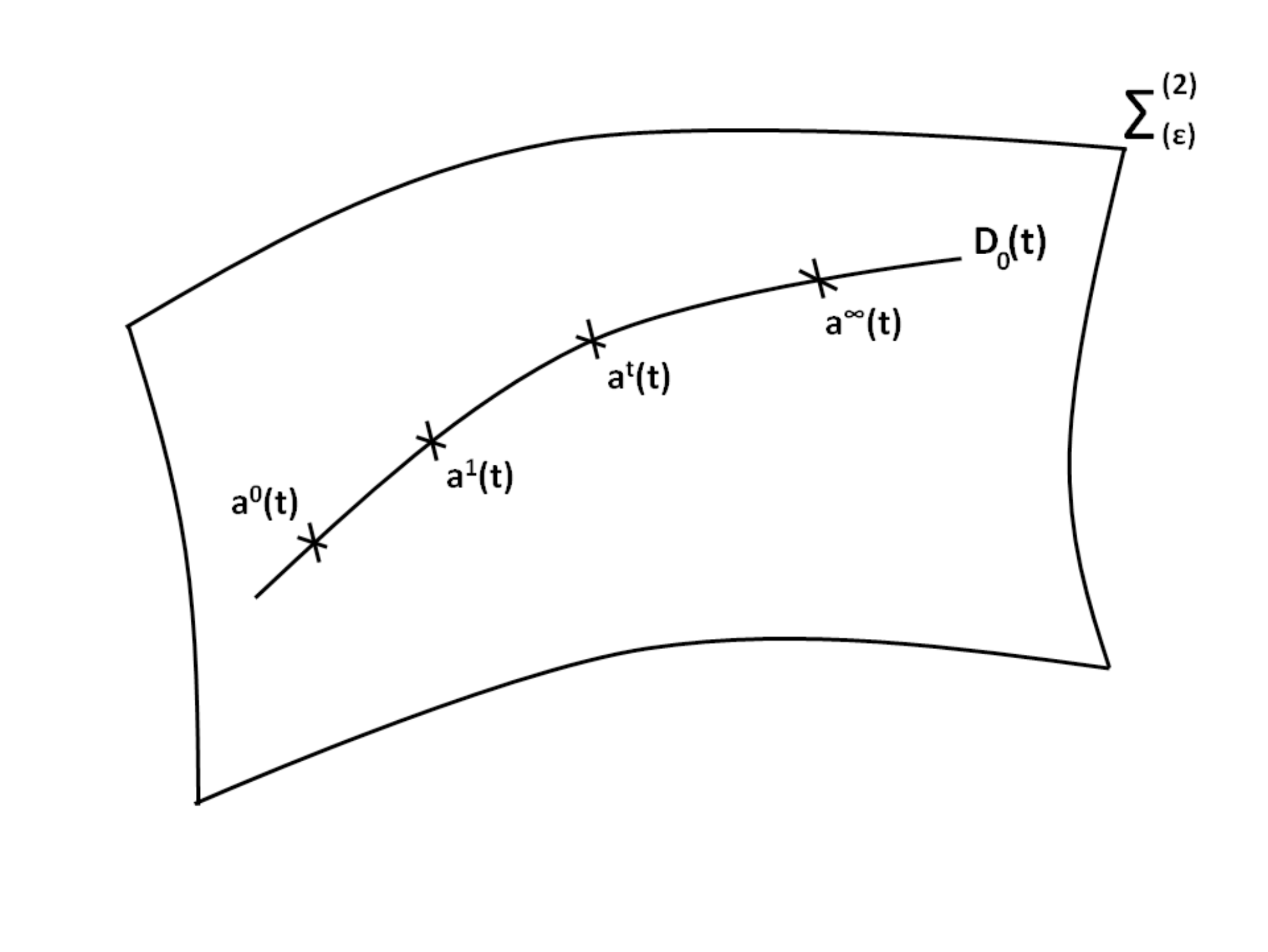}
\end{center}
\caption{The divisor $D_0(t)$ in the chart $W_2$.} \label{fig1}
\end{figure}
The holonomy group of the Painlev\'{e} VI foliation at infinity is, by definition,
the  holonomy group along the vertical leaf  $D_0=D_0(c)$, $c\neq 0,1, \infty$,
shown on fig.\ref{fig2}. The leaf $D_0(c)$ is a four-punctured Riemann sphere, the punctures corresponding to $a^0, a^1, a^c, a^\infty$. 
The holonomy group of the Painlev\'e VI foliation along $D_0(c)$ is   generated by three germs of analytic diffeomorphisms
\begin{equation}
\label{diffeos}
h_0, h_1, h_c
\end{equation}
corresponding to loops on $D_0(c)$ making one turn around $a^0, a^1, a^c$ respectively.

The main result of the paper, Theorem
\ref{main}, follows from Lemma \ref{mainlemma} and Lemma \ref{mainlemma2}
formulated bellow.

Let $G$ be a group with three generators $a,b,c$ and the following defining
relations
$$
a^2=b^2=c^2=(abc)^2=1 .
$$
An element $g\in G$ can be therefore represented by a word formed by the
letters $a,b,c$. The length $l(g)$ of a word $g\in G$ is the number of letters
in $g$, and only the equivalence class of $l(g)$ in
$\mathbb{Z}_2=\mathbb{Z}/2\mathbb{Z}$ is well defined. We get a homomorphism
$$
G\rightarrow \mathbb{Z}_2 : g \mapsto l(g)
$$
and let $G^0$ be its kernel.
\begin{lemma}
\label{mainlemma} The group $G^0$ is commutative.
\end{lemma}
\textbf{Proof.} The above lemma has a transparent geometric meaning : $G^0$ is
identified to the fundamental group of an elliptic curve which, as well known,
is commutative. To see this, consider the compact Riemann surface $S$ with
affine model
$$
S=\{(x,y): y^2=x(x-1)(x-t) \}, t\neq 0,1
$$
as well the natural projection
$$
\pi : S \rightarrow \mathbb{P} : (x,y) \rightarrow x .
$$
Let $\tilde{S}=S\setminus\{(x,0): x= 0,1,t,\infty\}$ and $\tilde{\mathbb{P}} =
\mathbb{P}\setminus \{0,1,t,\infty \}$. Let $P\in \tilde{S}$ and, by abuse of
notation,  $p=\pi(P)$. The fundamental group $\pi_1(\tilde{\mathbb{P}} ,p)$ is
the free group generated by $a,b,c$, where  $a,b,c$ are represented by closed
loops making one turn around $0,1,t$. With orientations appropriately chosen,
$abc$ is represented by a loop around $\infty$. The fundamental group
$\pi_1(\tilde{S} ,p)$ is  free with five generators. The projection
$$
\pi : \tilde{S} \rightarrow {\mathbb{P}} : (x,y) \rightarrow x .
$$
is a two sheeted covering which induces a monomorphism
$$\pi_* : \pi_1(\tilde{S} ,P) \rightarrow \pi_1(\tilde{\mathbb{P}} ,p)$$
such that $\pi_1(\tilde{\mathbb{P}} ,p) / \pi_* ( \pi_1(\tilde{S} ,P)) =
\mathbb{Z}_2$. This can be resumed in the following exact sequence of
homomorphisms
$$
1\rightarrow \pi_1(\tilde{S} ,P) \rightarrow \pi_1(\tilde{\mathbb{P}}
,p)\rightarrow \mathbb{Z}_2 \rightarrow 1 .
$$
 An element of the fundamental group $\pi_1(\tilde{S},P)$
represented by a closed loop which makes one turn around one of the
ramification points on $S$ is mapped by $\pi$ to $a^2$,$b^2$, $c^2$ or
$(abc)^2$. It follows that the induced homomorphism
$$
\pi_1(S,P) \rightarrow  G
$$
is well defined, where $G$ is the group defined above. The image of
$\pi_1(S,P)$ in $G$ consistes of words of even length and each word of even
length has a unique pre-image (lift of a closed loop with respect to the
projection). Therefore the following  sequence of homomorphisms is exact
$$
1\rightarrow \pi_1(S ,P) \rightarrow G \rightarrow \mathbb{Z}_2 \rightarrow 1
$$
and $Im (\pi_1(S ,P))=G^0$ is commutative.
 The lemma is proved.$\Box$

Next we apply the above Lemma to the holonomy group along $D_0(c)$. This group has three generators
 $h_0,h_1, h_c$, see (\ref{diffeos}), and let $h_\infty$ be the holonomy map associated to a closed loop making one turn around $\lambda=\infty$. If the orientations of the underlying closed loops are appropriately chosen, then
 $h_\infty = h_0h_1h_c$. 
\begin{lemma}
\label{mainlemma2}
$$
h_0^2=h_1^2=h_c^2=(h_{0}h_{1}h_{c})^2 = id .
$$
\end{lemma}
\textbf{Proof.} Recall that, according to section \ref{compactification} the
Okamoto surface $\bar{E}$ is obtained from the Hirzebruch surface $
\Sigma^{(2)}_{(\varepsilon)}$ after $ 8=4\times 2$ blow up's at the four
singular points $a^0, a^1, a^t, a^\infty$. It follows that the holonomy group
along the leaf $D_0=D_0(t)\subset \bar{E}$ coincides with the holonomy group of
the divisor leaf $D_0=D_0(t)\subset  \Sigma^{(2)}_{(\varepsilon)}$, see
fig.\ref{fig1} and fig.\ref{fig2}. Further, because of the symmetry of Painlev\'{e}
VI \cite[Proposition 2.2]{okamoto79} it suffices to show that $h_0^2=id$.

The idea of the proof of Lemma \ref{mainlemma2} is as follows. By Theorem
\ref{mainokamoto} consider the family of solutions $\lambda=\lambda(t),
\mu=\mu(t)$, with initial conditions
$$
(\lambda(t_0),\mu(t_0))\in D_2^0(t_0) .
$$
Here $\lambda, \mu$ are coordinates in appropriate chart on the variety
$\bar{E}_c$. Clearly these solutions are analytic for $t\sim t_0$. Upon a
successive contraction of the divisors $D_2^0$ and $D_1^0$ we get an infinite
family of analytic solutions $\lambda=\lambda(t), \mu=\mu(t)$ which tend to the
point $a^0$. The corresponding leaves are in fact holomorphic curves at $a^0$
\cite[Lemme $2_{VI}$]{okamoto79} which are therefore  parameterized by suitable convergent Puiseaux series
$$
\lambda \rightarrow (t(\lambda), \mu(\lambda)).
$$ 
The monodromy of these series when $\lambda$ makes one turn around the origin is readily computed to be an involution.
We claim that all leaves
"sufficiently close" to $D_0(c)$ are obtained in such a way, with some $t_0\sim
c$. If true, this would imply that $h_0^2=id$ as this holds true for the
monodromy map of the holomorphic curves through $a_0$.

 To make these considerations rigorous, consider the chart $W_2$ on the Hirzebruch surface, with coordinates $(\lambda_2,\mu_2)$,
see section \ref{compactification}. 
The Painlev\'{e} foliation along the leaf $D_0(c)=\{t=c,\mu_2=0\}$ on
$\bar{E}$, after $8=4\times2$ blow downs, is defined by (\ref{2prime}).
Until the end of this section we replace, as in section \ref{higher}  $(\lambda_2,\mu_2)$ by $(\lambda,\mu)$.
\begin{figure}
\begin{center}
\includegraphics[width=12cm]{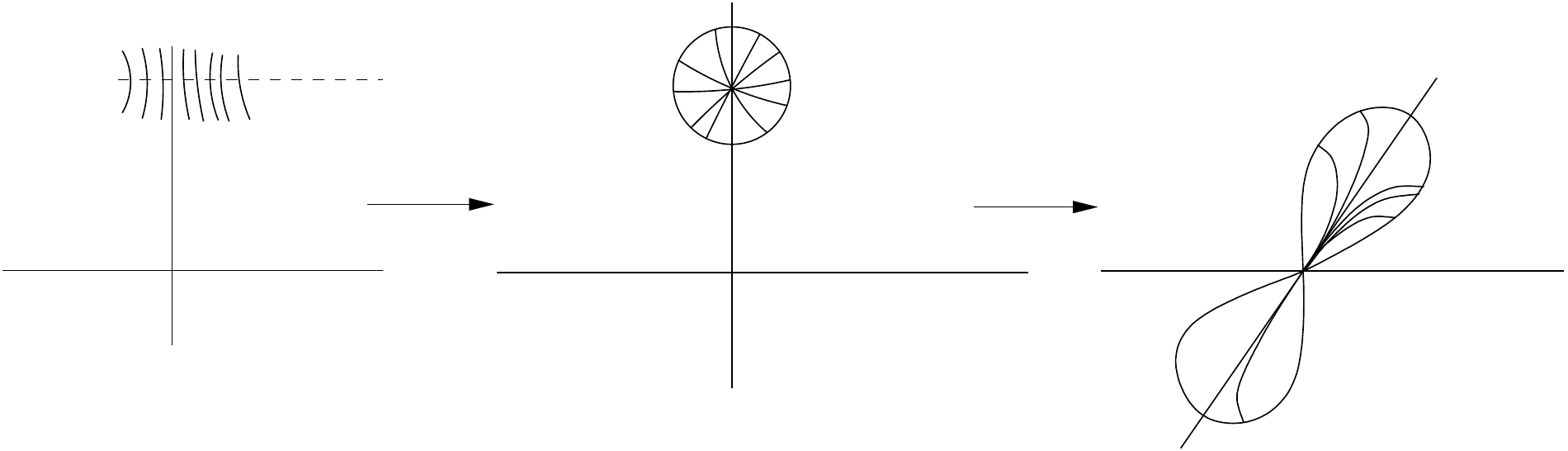}
\end{center}
\caption{Blowing down the divisors $D_2^0(c)$ and $D_1^0(c)$, and the domain
$U_c$} \label{fig3}
\end{figure}

Consider an open neighborhood $\tilde{U}_c$ of the divisor
$D_2^0=D_2^0(c)=\mathbb{P}^1$ shown on fig.\ref{fig2}. in the three-dimensional space $\bar{E}$. As the
Painlev\'{e} foliation is transversal to $D_2^0(c)$, then we shall suppose that
$\tilde{U}_c$ intersects any leaf of the foliation into an open disc, and that
$\tilde{U}_c$ is a union of such discs. After a contraction of $D_2^0(t)$ and
$D_1^0(t)$, $\forall t$, the neighborhood $\tilde{U}_c$  is transformed to a
cone-like domain  $U_c \subset \Sigma^{(2)}_{(\epsilon)}\times B . $ The
processus of blowing down the divisors $D_2^0(t)$ and $D_1^0(t)$, and the
effect on the domain $U_c\cap \{t=const.\}$ is shown on fig.\ref{fig3}.  
As the neighborhood $\tilde{U}_c$ is a union of regular leaves, then the contracted domain  $U_c$ is an union of leaves of (\ref{2prime}) intersecting at $a^0(t)$.
Each leaf of (\ref{2prime}) is therefore  a holomorphic curve at $a^0(t)$, tangent to the plane $\{ \lambda
+\kappa_0\mu=0 \}$ there. 
Another important feature of the contracted neighborhood $U_c$ is that it contains the  domain
\begin{equation}\label{cone}
\{ |\frac{\lambda}{\mu}+ \kappa_0|<\varepsilon, |t-c|<\delta ,
|\lambda|<\varepsilon, |\mu|<\varepsilon \}
\end{equation}
for all sufficiently small $\varepsilon,\delta >0$. The domain (\ref{cone}) is a direct product of the cone
$$
\{ (\lambda,\mu)\in \C^2 :  |\frac{\lambda}{\mu}+ \kappa_0|<\varepsilon, |\lambda|<\varepsilon, |\mu|<\varepsilon \}
$$
with vertex $(0,0,c)$ and axis 
$\lambda
+\kappa_0\mu=0$, and the disc $\{t\in \C : |t-c|<\delta \}$. 
Consider the following projection map
\begin{equation}
\label{pi}
\pi :  (\lambda,\mu,t) \mapsto (\lambda
+\kappa_0\mu,0,c)
\end{equation}
 from a neighborhood of the point $a^0(c)=(0,0,c)$ to the line $D_0(c)$. The pre-image of $\pi^{-1}(a^0(c))$ is the complex  two-plane 
 $ \{\lambda
+\kappa_0\mu=0\}$.

Consider a cross-section $\sigma=(\mathbb{C}^2,0)$ to $D_0(c)$, contained in
the plane  $\lambda=\varepsilon>0$. It is a (germ of a) complex two-dimensional disc centered at the origin on which the holonomy map is defined 
$$
h_0: \sigma \rightarrow \sigma .
$$
To define geometrically $h_0$, consider
  the path $\gamma=[0,\varepsilon]\subset
D_0(c)$ (a real interval) connecting $\lambda=\varepsilon$ to the origin on the
$\lambda$-plane $D_0(c)$. We claim that for any initial condition on $\sigma$
sufficiently close to $(\varepsilon,0,c)$ there is a lift $\Gamma$ of $\gamma$
along $\pi$, to a path contained in a leaf of the Painlev\'{e} foliation, starting
at the above initial point. Moreover, we claim that when $\lambda$ tends to
zero along $\gamma$, then the corresponding point of  $\Gamma$ tends to
$(0,0,t)$ for some $t$, $|t-c|\leq \varepsilon$. For this purpose we prove
first that $\gamma$ can be lifted at least until it intersects the cone-like
domain $U_c$. As the leaves in $U_c$ are  curves holomorphic at $a^0(c)$, then the result will
follow.

To lift $\gamma$ until it intersects $ U_c$ we construct a suitable compact
set, in which $\gamma$ can be lifted. Namely, let $K$ be the closure of the
following set
$$
 \{(\lambda,\mu,t): |\frac{\lambda}{\mu}+ \kappa_0|\geq
\delta, |t-c|\leq\delta,  |\lambda  + \kappa_0  \mu| \leq \delta \} .
$$
It is easily seen that $K$ is compact. The foliation (\ref{2prime}) is
transverse to the fibers of the map $\pi$ at a point $(\lambda,\mu,t)$ if
$$
 2E(t,\lambda)+F(t,\lambda)\mu + \kappa_0
\left[E_{\lambda}(t,\lambda)+F_{\lambda}(t,\lambda)\mu+G\mu^2\right]\mu \neq 0
.
$$
In a suitable neighborhood of the point $(0,0,c)$ we have
$$
2E(t,\lambda)+F(t,\lambda)\mu + \kappa_0
\left[E_{\lambda}(t,\lambda)+F_{\lambda}(t,\lambda)\mu+G\mu^2\right]\mu =
2t(\lambda+ \kappa_0 \mu)+ \dots
$$
where the dots stand for $O(|\lambda|^2+|\mu|^2)$ uniformly in $t\sim c$. We
conclude that when $(\lambda,\mu,t)\in K$ and belongs to a suitable
neighborhood of the point $(0,0,c)$, then the foliation (\ref{2prime}) is
transverse to the fibers of the map $\pi$, and hence
 the path $\gamma$ can be
lifted until its lift reaches the border of $K$. In the case when
$|t-c|<\delta$, this means that the path can be lifted until the cone-like
domain $U_c$ which is filled up by holomorphic curves (leaves of the
foliation). Therefore the path  can be further lifted until the origin and the
claim is proved.

It remains to show that in the course of the lifting $|t-c|<\delta$ holds true.
This follows after integrating the differential
$$
\frac{dt}{t(t-1)} = \dfrac{1}{2E(t,\lambda)/\mu +F(t,\lambda)} d\lambda 
$$
along the path $\gamma$. Indeed, when $(\lambda,\mu,t)\in K$ and belongs to a
suitable neighborhood of the point $(0,0,c)$, then
$$
\frac{2E(t,\lambda)}{\mu}\sim -2\kappa_0 c, F(t,\lambda)  \sim \kappa_0 c
$$
and hence
$$
\frac{2E(t,\lambda)}{\mu} + F(t,\lambda) \sim -\kappa_0 c
$$
is bounded   from zero, provided that $\kappa_0\neq 0$. 

We conclude that every path can be lifted until it crosses the domain $U_c$ in which the
leaves of the foliation are holomorphic curves. The monodromy of the Puiseaux series
$t=t(\lambda'), \mu=\mu(\lambda')$, when $\lambda'=\lambda +\kappa_0\mu$ makes one turn around the
origin is easily described : it is an involution. This  follows geometrically
from the fact, that the divisor $D^0_2(t)$ is obtained after two blow up's from
$a^0(t)$. Analytically, this means that if $z$ is a local coordinate on the projective line $D_2^0(t)$, then $U_c$ is an union of holomorphic curves (leaves) parameterized by $z$ and $t$, such that for  fixed $z,t$ we have
$$
\lambda' = \lambda +\kappa_0\mu= z \mu^2 +O(\mu^3).
$$
where $-\kappa_0$ is the coordinate of $D_1^0(t) \cap D_2^0(t)$ and $z \in D_2^0(t)$ is the intersection point of the leaf and $D_2^0(t)$.
Therefore $\mu$ is an analytic function in $\sqrt{\lambda'}$ and the result follows.
Finally, we note that the holonomy map $h_0$ depends analytically on the parameters of the Painlev\'{e} foliation. As $h_0^2=id$ for $\kappa_0\neq0$ then this holds true for all $\kappa_0$. To
resume, we proved

\begin{proposition}
There exists  a neighborhood of the point $a^0(c)=(0,0,c)$, such that every leaf of
the Painlev\'{e} foliation is a holomorphic curve, which is a ramified two-sheeted
covering of the divisor $D_0(c)=\{\mu=0\}$ along the projection map $\pi$
(\ref{pi}), with ramification  point $(\lambda=\mu=0)$.
\end{proposition}

The above Proposition generalizes \cite[Okamoto, Lemma $2_{VI}$]{okamoto79} and
implies Lemma \ref{mainlemma}.$\Box$

\textbf{Proof of Theorem \ref{main}.} According to Lemma \ref{mainlemma2} each
element of the holonomy group is a word made with the letters $h_0, h_1, h_c$.
By Lemma \ref{mainlemma} the subgroup of the holonomy group formed by words of
even length is commutative.

\vspace{2ex} \noindent {\bf Acknowledgment.} Part of this paper was written while the first author was visiting the Paul Sabatier University of Toulouse. He is obliged for the hospitality.


\begin{thebibliography}{00}


\bibitem{calo09} 
S. Cantat,  F. Loray,  Holomorphic dynamics, Painlev\'e VI equation, and character varieties ; Annales de l'Institut Fourier, Vol. 59 no. 7 (2009), p. 2927-2978.
\bibitem{casa09} G. Casale,
Une preuve galoisienne de l'irr\'eductibilit\'e au sens de
	Nishioka-Umemura de la 1ere \'equation de Painlev\'e
	Ast\'erisque 323 (2009)
	\bibitem{casa07} G. Casale,
		The Galois groupoid of Picard-Painlev\'e sixth equation
	RIMS Kokyuroku Bessatsu vol B2 (2007)
	\bibitem{casa12a} G. Casale,
		An introduction to Malgrange pseudogroup
	To appear in S\'eminaires et Congres (2012)
	\bibitem{casa12b}	G. Casale,
	Liouvillian first integrals of differential equations
	To appear in Proceedings of the Banach Center (2012)
	


\bibitem{hain87} R.M. Hain, The geometry of the mixed Hodge structure on the fundamental group, \textit{Proceedings of Symposia in Pure Mathematics} {\bf 46} (1987) 247-281.

\bibitem{horozov07} E. Horozov, T. Stoyanova, Non-Integrability of some Painlev\'e VI equations and dilogarithms, \textit{Regular \& Chaotic Dynamics} {\bf 12} (6) (2007) 622-629.

\bibitem{8} K. Iwasaki, H. Kimura, S. Shimomura, M. Yoshida, From Gauss to Painlev\'e, A Modern Theory of Special Functions, \textit{Aspects of Mathematics}, E16, Vieweg \& Shon, 1991.

\bibitem{malgrange} B. Malgrange, Le  groupoide de Galois d'un feuilletage, Monographie 38 vol 2, L'Enseignement Math\'ematique (2001).
\bibitem{mora01}
Juan~J. Morales-Ruiz and Jean~Pierre Ramis.
\newblock Galoisian obstructions to integrability of {H}amiltonian systems.
  {I}, {II}.
\newblock {\em Methods Appl. Anal.}, 8(1):33--95, 97--111, 2001.


\bibitem{ramis07} J.J. Morales-Ruiz, J.P. Ramis, C. Simo, Integrability of Hamiltonian systems and differential Galois groups of higher variational equations,
\textit{Ann. Scient. \'Ec. Norm. Sup.} {\bf 40} (4) (2007) 845-884.

\bibitem{morales06} J.J. Morales-Ruiz, A remark about the Painlev\'e transcendents, \textit{S\'eminaires \& Congr\'es} {\bf 14} (2006) 229-235.

\bibitem{morales99} J.J. Morales-Ruiz, Differential Galois theory and non-integrability of Hamiltonian systems, \textit{Progress in Mathematics}, {\bf 179}, Birkhauser, 1999.
\bibitem{mora10}
Morales-Ruiz, Juan J.; Ramis, Jean-Pierre Integrability of dynamical systems through differential Galois theory: a practical guide, in \textit{Differential algebra, complex analysis and orthogonal polynomials}, 143 - 220, Contemp. Math., 509, Amer. Math. Soc., Providence, RI, 2010
\bibitem{okamoto79} K. Okamoto, Sur les feuilletages associ\'{e}s aux \'{e}quations du second ordre \`{a} points critiques  fixes de P. Painlev\'{e}, \textit{Japan J. Math.} \textbf{5} (1979) 1-79.
\bibitem{stoy09}
Tsvetana Stoyanova,
\newblock Non-integrability of {P}ainlev\'e {VI} equations in the {L}iouville
  sense.
\newblock {\em Nonlinearity}, 22(9):2201--2230, 2009.
\bibitem{tret79}
Carol Tretkoff,  Marvin  Tretkoff, 
Solution of the inverse problem of differential Galois theory in the classical case.
Amer. J. Math. 101 (1979), no. 6, 1327-1332. 
\bibitem{wata98}
Humihiko Watanabe, Birational canonical transformations and classical solutions of
the sixth Painlev\'e equation, 
\textit{Ann. Scuola Norm. Sup. Pisa Cl. Sci.} 
(4), 27 (1998), no. 3-4, 379-425 (1999).
\end{thebibliography}
\end{document}